\documentclass[12pt]{article}
\usepackage[margin=0.8in]{geometry}
\usepackage{amsmath,amssymb,amsthm}
\usepackage{color}
\usepackage{parskip}
\usepackage{setspace}
\usepackage{etoolbox}
\patchcmd{\thebibliography}{*}{}{}{}

\newtheorem{theorem}{Theorem}[section]

\newtheorem{lemma}[theorem]{Lemma}

\newtheorem{definition}[theorem]{Definition}

\newtheorem{remark}[theorem]{Remark}

\newtheorem*{definition*}{Definition}

\usepackage{amssymb,latexsym,amsmath,epsfig,amssymb,amsthm, amsmath,color}

\begin{document}
\title{A note on conditional expanders over prime fields}
\author{Mozhgan Mirzaei\thanks{Department of Mathematics,  University of California at San Diego, La Jolla, CA, 92093 USA.  Supported by NSF grant DMS-1800746. Email:
{\tt momirzae@ucsd.edu}.}}
\maketitle


\date{}
\begin{abstract}
Let $\mathbb{F}_p$ be a prime field of order $p,$ and $A$ be a set in $\mathbb{F}_p$ with $|A| \leq p^{1/2}.$ In this note, we show that
\[\max\{|A+A|, |f(A, A)|\}\gtrsim |A|^{\frac{6}{5}+\frac{4}{305}},\]
where $f(x, y)$ is a \textit{non-degenerate} quadratic polynomial in $\mathbb{F}_p[x, y].$ This improves a recent result given by Koh, Mojarrad, Pham, Valculescu (2018).
\end{abstract}

\section{Introduction}
Let $A$ be a set of integers. The sum and product sets are defined as follows:
\begin{align*}
A+A =& \{a+b\colon a, b \in A\}\\
A\cdot A =&\{ab\colon a,b \in A\}.
\end{align*}

Throughout this paper, by $X\gg Y,$ we mean $X\ge C_1Y$ for some positive constant $C_1$, and $X\sim Y$ means that $X\gg Y$ and $Y\gg X$, by $X \gtrsim Y$ we mean $X\gg (\log Y)^{-C_2} Y$ for some positive constant $C_2$.

Erd\H{o}s and Szemer\'edi \cite{es} proved that for any finite set $A \subset \mathbb Z$, we have
$$\max\{ |A+A|, |A \cdot A|\} \gg |A|^{1+\varepsilon}$$
for some positive constant $\varepsilon$.  In the setting of finite fields, a similar result has been derived by Bourgain, Katz, and Tao \cite{bourgain-katz-tao}. They showed that for any set $A \subset \mathbb F_{p}$, where $p$ is a prime and $p^{\delta} < |A| < p^{1-\delta}$ for some $\delta >0,$ one has
$$\max\{ |A+A|, |A \cdot A|\} \geq C_{\delta}|A|^{1+\varepsilon},$$ for some $\varepsilon=\varepsilon(\delta) > 0$. 
We note here that in the result of Bourgain, Katz, and Tao \cite{bourgain-katz-tao}, it is difficult to determine the relation between $\varepsilon$ and $\delta$. 

Hart, Iosevich, and Solymosi \cite{his} developed Fourier analysis tools to obtain a bound over arbitrary finite fields that gives an explicit dependence of $\varepsilon$ on $\delta$ as follows. 
\bigskip
\begin{theorem}[\textbf{Hart-Iosevich-Solymosi,} \cite{his}]\label{hiss}
Let $ \mathbb{F}_q$ be an arbitrary finite field of order $q$, and let $A \subset \mathbb{F}_q$. Suppose $|A + A| = m$ and $|A \cdot A | = n$, then we have
\begin{equation}\label{eq:his}
|A|^3 \leq \frac{ c m^2 n |A| }{ q} + c q^{1/2} mn,
\end{equation}
for some positive constant $c$. 
\end{theorem}
By a direct computation, Theorem \ref{hiss} is non-trivial when $|A|\gg q^{1/2}.$ For $|A|\sim q^{7/10},$ we have 
\[\max\left\lbrace |A+A|, |A\cdot A|\right\rbrace \gg  |A|^{8/7}.\]

Using exponential sums, 
Garaev \cite{garaev} obtained the following improvement.
\bigskip
\begin{theorem}[\textbf{Garaev}, \cite{garaev}]\label{vinh cor}
Let $\mathbb{F}_q$ be a finite field of order $q$ and $A$ be a set in $\mathbb{F}_q$. 
\begin{enumerate}
\item If $q^{1/2} \ll |A| \ll q^{2/3}$, then 
\[
\max \left\lbrace |A + A | , |A \cdot A | \right\rbrace \gg \frac{ |A|^2 }{q^{1/2}} .
\]
\item If $|A| \gg q^{2/3}$, then 
\[
\max \left\lbrace |A + A | , |A \cdot A | \right\rbrace \gg ( q |A| )^{1/2} .\]
\end{enumerate}
\end{theorem}
\bigskip
Hence, if $|A|=q^{\alpha}$, then we have 
\[\max\left\lbrace |A+A|, |A\cdot A|\right\rbrace\gg |A|^{1+\alpha'}, \]
where $\alpha'=\min\left\lbrace 1-1/2\alpha, (1/\alpha-1)/2\right\rbrace$. If $\alpha$ is very small, say $\alpha\le 18/35$, and $q$ is a prime number, then Rudnev, Shakan,  and Shkredov \cite{sha} proved the following. 
\bigskip
\begin{theorem}[\textbf{Rudnev-Shakan-Shkredov}, \cite{sha}]\label{thm20}
Let $\mathbb{F}_p$ be a prime field of order $p$. Let $A$ be a set in $\mathbb{F}_p$. Suppose that $|A|\ll p^{\frac{18}{35}}$, then we have 
\[\max\{|A+A|, |A\cdot A|\}\gg |A|^{1+\frac{2}{9}-o(1)}.\]
\end{theorem} 
\bigskip
This theorem improves the earlier exponents $39/32$ due to  Chen, Kerr, and Mohammadi \cite{chen} and $6/5$ due to Roche-Newton, Rudnev, and Shkredov \cite{RRS}. 
\bigskip
 \begin{definition} \label{def:deg} Let $\mathbb{F}_p$ be a prime field.  A polynomial $f(x, y) \in \mathbb{F}_p[x, y]$ is {\it degenerate}
 if it is  of the form
$Q(L(x, y))$ where $Q$ is an one-variable polynomial and $L$ is
a linear form in $x$ and  $y$. \end{definition}

In a recent work, a very general bound for quadratic polynomials has been given by Koh, Mojarrad, Pham, and Valculescu \cite{koh1}. 

\bigskip
\begin{theorem}[\textbf{Koh-Mojarrad-Pham-Valculescu},  \cite{koh1}]\label{kohj}
Let $\mathbb{F}_p$ be a prime field of order $p,$ and let $f(x, y)\in \mathbb{F}_p[x,y]$ be a non-degenerate quadratic polynomial. For $A\subset\mathbb{F}_p$ with $|A|\ll p^{5/8}$, we have 
\[\max\{|A+A|, |f(A, A)|\}\gg |A|^{6/5}.\]
\end{theorem}

In this paper, we use methods in \cite{sha, ss} to improve this theorem as follows. 
\bigskip
\begin{theorem}\label{thm1}Let $f(x, y)\in \mathbb{F}_p[x, y]$ be a  non-degenerate quadratic polynomial. For $A\subset\mathbb{F}_p$ with $|A|\ll p^{1/2}$ we have
\[\max\{|A+A|, |f(A, A)|\}\gtrsim |A|^{\frac{6}{5}+\frac{4}{305}}.\]
\end{theorem}
\section{Proof of Theorem \ref{thm1}}
For $A, B\subset\mathbb{F}_p$, let $\mathcal{E}_4^+(A, B)$ be the number of tuples $(a_1, a_2, a_3, a_4, b_1, b_2, b_3, b_4)\in A^4\times B^4$ such that 
\[a_1-b_1=a_2-b_2=a_3-b_3=a_4-b_4.\]

For $A\subset \mathbb{F}_p$, we define 
\[d_4^+(A):=\sup_{B\ne \emptyset} \frac{\mathcal{E}_4^+(A, B)}{|A||B|^{3}}.\]

Note that since $\mathcal{E}_4^+(A, B)\ge \mathcal{E}_4^+(A, A)\ge |A|^4$, we have $d_4^+(A)\ge 1$. 

It has been observed in \cite{ss} that the $\sup$ is taken over all sets $B$ with $|B|\le |A|^{3/2}.$ Indeed, if $|B| \geq |A|^{3/2}$, then \[d_4^+(A)=sup_{B\ne \emptyset} \frac{\mathcal{E}_4^+(A, B)}{|A||B|^{3}}\leq \frac{|A|^4|B|}{|A||B|^3} \leq 1,\] 
a contradiction. 

In \cite{ss}, Shakan and Shkredov proved that 
\[d_4^+(A)\ll \frac{|A\cdot A|^2}{|A|^2}\]
whenever $|A|\le p^{3/5}$. They also derived the following lemma. 
\bigskip
\begin{lemma}[\cite{ss}]\label{hai}
For $A\subset \mathbb{F}_p$, we have 
\[d_4^+(A)\gtrsim \frac{|A|^{48/13}}{|A+A|^{35/13}}.\]
\end{lemma}
\bigskip
In this paper, we will give an upper bound of $d_4^+(A)$ in terms of $|f(A, A)|$ for any non-degenerate quadratic polynomial $f$ as follows. 
\bigskip

\begin{lemma}\label{lm1}
Let $f(x, y)\in \mathbb{F}_p[x, y]$ be a  non-degenerate quadratic polynomial. For $A\subset \mathbb{F}_p$ with $|A|\ll p^{1/2}$, we have 
\[d_4^+(A)\ll \frac{|f(A, A)|^{2}}{|A|^{2}}.\]
\end{lemma}
\bigskip
To prove lemma \ref{lm1}, we use the following result in \cite{KMPS}.
\bigskip
\begin{theorem}[\cite{KMPS}]\label{thm12}
Let $f\in \mathbb{F}_p[x,y,z]$ be a quadratic polynomial that depends on each variable and that does not have the form $g(h(x)+k(y)+l(z))$. For $A, B, C\subset \mathbb{F}_p$ with $|A||B||C|\ll p^2$, let $E$ be the number of tuples $(a, b, c, a', b', c')\in (A\times B\times C)^2$ such that $f(a, b, c)=f(a', b', c')$. Then we have 
\[E\ll (|A||B||C|)^{3/2}+(|A|+|B|+|C|)(|A||B||C|)+|B|^2|C|^2.\]
\end{theorem}

\paragraph{Proof of Lemma \ref{lm1}:}
Let $B\subset \mathbb{F}_p$ be a set maximizing $d_4^+(A).$
By a dyadic decomposition, there exist a number $t>0$ and a set $D_t:=\{x\colon r_{A-B}(x)\ge t\}$ such that 
\[\mathcal{E}_4^+(A, B)\sim |D_t|t^4.\] 
Without loss of generality, we assume that  $f(x, y)=ax^2+by^2+cxy+dx+ey$ with $a\ne 0$. Let $f'(u, v, w):=f(u+v, w)$. 

Since $f(x, y)$ is a non-degenerate polynomial, by an elementary calculation, we have $f'(x, y, z)$ is not of the form $g(h(x)+k(y)+l(z))$ for some polynomials $g, h, k, l$.  

Consider the following equation
\begin{equation}\label{eqmay}f'(u, v, w)=t, \end{equation}
with $u\in D_t, v\in B, w\in A, t\in f(A, A)$. 

It easy to check that the number of solutions of the equation (\ref{eqmay}) is at least $|D_t|t|A|$. Now by the Cauchy-Schwarz inequality, we have 
\begin{equation}\label{e1x}|D_t|t|A|\ll |f(A, A)|^{1/2}E^{1/2},\end{equation}
where $E$ is the number of tuples $(u, v, w, u', v', w')\in (D_t\times B\times A)^2$ such that 
\[f(u, v, w)=f(u', v', w').\]
Suppose $|D_t||A||B|\ll p^2$. We now consider the following cases:

{\bf Case $1$:} If $|D_t|\le |B|$, then Theorem \ref{thm12} with $A:=D_t, B:=B, C:=A$ gives
\[E\ll (|D_t||B||A|)^{3/2}+(|D_t|+|B|+|A|)(|D_t||B||A|)+|B|^2|A|^2.\]

{\bf Case $2$:} If $|D_t|\ge |B|$, then Theorem \ref{thm12} with $A:=B, B:=D_t, C:=A$ gives 
\[E\ll (|B||D_t||A|)^{3/2}+(|B|+|D_t|+|A|)(|B||D_t||A|)+|D_t|^2|A|^2.\]

These cases can be handled in the same way. Therefore, without loss of generality, we assume that we are in the first case, i.e. $|D_t|\le |B|$. 

We also can assume that $|D_t||A|\ge |B|$, otherwise, using the fact that $t\le |A|, |B|$, we have 
\[\frac{|D_t|t^4}{|A||B|^3}\ll \frac{|B|t^3}{|A||B|^3}\ll 1\le \frac{|f(A, A)|^2}{|A|^2}.\]
Similarly, we assume that $|D_t||B|\ge |A|$. Furthermore, if $(|A||B|)^3\le |D_t|^5$, then we have $|A|^3|B|^3\le |B|^5$. This implies that $|B|^2\ge |A|^3$. This is not possible since $|B|\le |A|^{3/2}$. Thus we can assume that $|A|^3|B|^3 \geq |D_t|^5$. With these assumptions, we obtain 
\[ E\ll (|D_t||A||B|)^{3/2}+(|B||A|)^2.\]
Without loss of generality, let us assume $(|D_t||A||B|)^{3/2} \geq (|B||A|)^2.$ (As otherwise, $|B||A| \leq |D_t|^3.$ Hence $\frac{|D_t|t^4}{|A||B|^3} \leq \frac{t^4}{|D_t|^2|B|^2} \leq \frac{t^4}{D^4_t} \leq 1 \leq \frac{|f(A,A)|^2}{|A|^2},$ and we are done.)
Therefore, 
\[|D_t|t^4\ll \frac{|f(A, A)|^2|B|^3}{|A|},\]
and we are done by the definition of $d_4^+(A)$. \\
Suppose $|D_t||A||B|\gg  p^2$, then one can apply a point-plane incidence bound due to Vinh \cite{vinh} to obtain an upper bound of $E$. More precisely, we have 
\[E\ll \frac{|D_t|^2|A|^2|B|^2}{p}.\]
Substituting this inequality to (\ref{e1x}) we get 
\begin{equation}\label{eq90}pt^2\le |f(A, A)||B|^2.\end{equation}
If $t^3\le \frac{|B|^2|f(A, A)|^2}{|A|^2}$, then we can use the fact that $|D_t|t\le |A||B|$ to obtain the desired bound for $d_4^+(A)$. Otherwise, we obtain 
\begin{equation}\label{eq91}t\ge \frac{|B|^{2/3}|f(A, A)|^{2/3}}{|A|^{2/3}}.\end{equation}
Therefore, it follows from (\ref{eq90}) and (\ref{eq91}) that 
\[p \frac{|B|^{4/3}|f(A, A)|^{4/3}}{|A|^{4/3}}\le pt^2\le |f(A, A)||B|^2.\]
Thus, 
\[|A|^6\ge |A|^3|B|^2\gg p^3,\]
which gives us $|A|\ge p^{1/2}$. Therefore, we have  proved that 
\[d_4^+(A)\ll \frac{|f(A, A)|^2}{|A|^2},\]
whenever $|A|\le p^{1/2}$.$\hfill \square$

Theorem \ref{thm1} follows by combining Lemma \ref{hai} and Lemma  \ref{lm1}.

\bigskip
\begin{remark}
It is clear that if $f(x, y)=xy$, then Theorem \ref{thm1} is weaker than Theorem \ref{thm20}. In our general setting, the main difficulty arises when we want to give an upper bound for $\mathcal{E}_2^+(A, A-A)$ in terms of $|f(A, A)|$, where $\mathcal{E}_2^+(A, B)$ is the number of tuples $(a_1, a_2, b_1, b_2)\in A^2\times B^2$ such that $a_1-b_1=a_2-b_2$.  For all non-degenerate quadratic polynomials, it seems very difficult to give such an upper bound, but for some special families of polynomials it is possible. For instance, if $f(x, y)=g(x)(h(x)+y)$ is a function defined on $\mathbb{F}_p^*\times \mathbb{F}_p^*$, where $g, h\colon \mathbb{F}_p^*\to \mathbb{F}_p^*$ are arbitrary functions, then one can follow the proof of \cite[Theorem $1.6$]{mo} to derive the following:
\[\mathcal{E}_2^+(B, C)\ll |A|^{-2}\left(|f(A, B)|^{3/2}|A|^{3/2}|C|^{3/2}+k|f(A, B)|A||C|\right),\]
where $k\le \max\{|A|, |C|, |f(A, B)|\}$ under the assumption $|f(A, B)||A||C|\ll p^2$. 

\end{remark}

\textbf{Acknowledgments.} I would like to thank Ben Lund for reading the manuscript carefully and for helpful comments.

\end{document}